\documentclass[12pt]{article}
\usepackage{mathrsfs}
\usepackage{epic,eepic,epsf,epsfig}
\usepackage{amsfonts,srcltx,mathrsfs}
\usepackage{multirow}
\usepackage{amsmath}
\usepackage{amssymb}
\usepackage{amsbsy}
\usepackage{graphicx}
\usepackage{amsfonts}%,srcltx}
\usepackage{color}
\usepackage{setspace}

\newtheorem{lem}{Lemma}[section]%
\newtheorem{theorem}[lem]{Theorem}%

\def\nd{\mathrel{\bigm|\kern-.7em/}}

\def\f{\noindent}

\def\P\GammaL{\hbox{\rm P\GammaL}}

\def\mod{\hbox{\rm mod }}

\begin{document}
\title{More results on the spectral radius of graphs with no odd wheels}

\footnotetext{E-mails: zhangwq@pku.edu.cn}

\author{Wenqian Zhang\\
{\small School of Mathematics and Statistics, Shandong University of Technology}\\
{\small Zibo, Shandong 255000, P.R. China}}
\date{}
\maketitle

\begin{abstract}
For a graph $G$, the spectral radius $\lambda_{1}(G)$ of $G$ is the largest eigenvalue of its adjacency matrix. An odd wheel $W_{2k+1}$ with $k\geq2$ is a graph obtained from a cycle of order $2k$ by adding a new vertex connecting to all the vertices of the cycle. Let ${\rm SPEX}(n,W_{2k+1})$ be the set of $W_{2k+1}$-free graphs of order $n$ with the maximum spectral radius. Very recently, Cioab\u{a}, Desai and Tait \cite{CDT2} characterized the graphs in ${\rm SPEX}(n,W_{2k+1})$ for sufficiently large $n$, where $k\geq2$ and $k\neq4,5$. And they left the case $k=4,5$ as a problem. In this paper, we settle this problem. Moreover, we completely characterize the graphs in ${\rm SPEX}(n,W_{2k+1})$ when $k\geq4$ is even and $n\equiv2~(\mod4)$ is sufficiently large. Consequently, the graphs in ${\rm SPEX}(n,W_{2k+1})$ are characterized completely for any $k\geq2$ and sufficiently large $n$.

\bigskip

\f {\bf Keywords:} spectral radius; spectral extremal graph; odd wheel; walk.

\end{abstract}

\baselineskip 17 pt

\section{Introduction}

All graphs considered in this paper are finite and undirected.
For a graph $G$, let $\overline{G}$ be its complement. The vertex set and edge set of $G$ are denoted by $V(G)$ and $E(G)$, respectively. For a vertex $u$, let $d_{G}(u)$ be its degree. Let $\Delta(G)$ be the maximum degree of $G$. For two disjoint subsets $S,T$ of $V(G)$, let $e_{G}(S,T)$ be the number of edges between $S$ and $T$ in $G$. Let $G[S]$ be the  subgraph of $G$ induced by $S$, and let $G-S=G[V(G)-S]$.  For two vertices $u$ and $v$, we say that $u$ is a neighbor of $v$ or $u\sim v$, if they are adjacent in $G$.   For $\ell\geq2$ vertex-disjoint graphs $G_{1},G_{2},...,G_{\ell}$, let $\cup_{1\leq i\leq \ell}G_{i}$ be the disjoint union of them. Let $G_{1}\vee G_{2}\vee\cdots\vee G_{\ell}$ be the join obtained from $\cup_{1\leq i\leq \ell}G_{i}$ by connecting each vertex in $G_{i}$ to each vertex in $G_{i+1}$ for any $1\leq i\leq \ell-1$.  For a certain integer $n$, let $K_{n},C_{n}$ and $M_{n}$ be the complete graph, the cycle and a perfect matching of order $n$, respectively. For any terminology used but not defined here, one may refer to \cite{CRS}.

 Let $G$ be a graph with vertices $v_{1},v_{2},...,v_{n}$. The {\em adjacency matrix} of $G$ is  $A(G)=(a_{ij})_{n\times n}$, where $a_{ij}=1$ if  $v_{i}\sim v_{j}$, and $a_{ij}=0$ otherwise. The {\em spectral radius} $\lambda_{1}(G)$ of $G$ is the largest eigenvalue of its adjacency matrix. By Perron--Frobenius theorem,  $\lambda_{1}(G)$ has a non-negative eigenvector (called Perron vector) and a positive eigenvector if $G$ is connected.  
For a graph $F$, a graph $G$ is call $F$-free if $G$ does not contain $F$  as a subgraph. Let ${\rm EX}(n,F)$ be the set of $F$-free graphs of order $n$ with the maximum number of edges, and let ${\rm SPEX}(n,F)$ be the set of $F$-free graphs of order $n$ with the maximum spectral radius.

 In 2010, Nikiforov \cite{N1} presented a spectral version of Tur\'{a}n-type
problem: what is the maximum spectral radius of an $F$-free graph of order $n$?
 In recent years, this problem has attracted many researchers (see \cite{BDT,CFTZ,CDT1,CDT2,CDT3,DKLNTW,FLSZ,FTZ,LL,LP1,LP2,LZZ,NWK,WKX,ZHL,ZL}). 
Let $F$ be a graph of chromatic number $r+1\geq3$. $F$ is called {\em color-critical}, if there is an edge $e$ of $F$ whose deletion induces a graph with chromatic number $r$. Simonovits \cite{S} showed that the only graph in ${\rm EX}(n,F)$ is the Tur\'{a}n
graph with $r$ parts when $n$ is sufficiently large. Nikiforov \cite{N2} proved that ${\rm SPEX}(n,F)={\rm EX}(n,F)$ for large $n$.
For $k\geq4$, let $W_{k}=K_{1}\vee C_{k-1}$, which is called a wheel of order $k$. Clearly, $W_{2k}$ for $k\geq2$ is a color-critical graph of chromatic number 4, and thus these theorems apply. For odd wheels $W_{2k+1}$ with $k\geq2$ and large $n$, the graphs in ${\rm EX}(n,W_{2k+1})$ were characterized by Dzido and  Jastrzebski  for $k=2$ \cite{DJ}, and by Yuan for $k\geq3$ \cite{Y}.
Recently, Zhao, Huang and Lin \cite{ZHL} determined the maximum spectral radius of graphs of order $n$ which forbid all wheels.
Very recently, Cioab\u{a}, Desai and Tait \cite{CDT2} characterized the graphs in ${\rm SPEX}(n,W_{2k+1})$ for sufficiently large $n$. For $k=2$, they prove the following  conclusion (see Theorem 1.2 in \cite{CDT2} and its proof).

\begin{theorem}{\rm(Cioab\u{a}, Desai and Tait \cite{CDT2})}
For sufficiently large $n$, the only graph in $\rm{SPEX}(n,W_{5})$ is obtained from the complete bipartite graph with parts $L$ and $R$ by embedding a maximum matching in each part of $L$ and $R$, where $|L|=\frac{n}{2}+1$ for $n\equiv2~(\mod4)$, and $|L|=\lceil\frac{n}{2}\rceil$ otherwise.
\end{theorem}

For $k\geq3$, a graph is called nearly $(k-1)$-regular if all the vertices have degree $k-1$ except one vertex with degree $k-2$. As in \cite{CDT2}, let $\mathcal{U}_{k,n}$ be the family of $(k-1)$-regular or nearly $(k-1)$-regular graphs
of order $n$ which do not contain a path of order $2k-1$. In fact, by Lemma \ref{bounded order} we see that $\mathcal{U}_{k,n}$ is exactly the family of $(k-1)$-regular or nearly $(k-1)$-regular graphs
of order $n$, in which each component has at most $2k-2$ vertices.

\begin{theorem}{\rm (Cioab\u{a}, Desai and Tait \cite{CDT2})} \label{CDT}
Let $k\geq3$ and $k\neq4,5$. For sufficiently large $n$, if $G\in\rm{SPEX}(n,W_{2k+1})$, then $G$ is obtained from
 a complete bipartite graph with parts $L$ and $R$ of size $\frac{n}{2}+s$ and $\frac{n}{2}-s$ respectively, by embedding a graph from
$\mathcal{U}_{k,\frac{n}{2}+s}$ in $G[L]$ and exactly one edge in $G[R]$. Furthermore, $|s|\leq1$.
\end{theorem}

At the end of \cite{CDT2}, the authors states that "We believe that the extremal graphs (in Theorem \ref{CDT}) when $k\in\left\{4,5\right\}$  have the same structure,
and it would be interesting to prove this".  Our first result is to confirm this.

\begin{theorem}\label{main1}
Theorem \ref{CDT} also holds for $k\in\left\{4,5\right\}$.
\end{theorem}

In the section 6 of \cite{CDT2}, the authors also determined the specified values of $|L|$ (or $s$) in Theorem \ref{CDT}. For odd $k$, $|L|=\lceil\frac{n}{2}\rceil$. When $k$ is even, $|L|$ is
constrained as follows.\\
 $|L|=\frac{n}{2}$ for $n\equiv0~(\mod4)$,\\
$|L|=\lfloor\frac{n}{2}\rfloor$ for $n\equiv1~(\mod4)$,\\
$|L|\in\left\{\frac{n}{2},\frac{n}{2}+1\right\}$ for $n\equiv2~(\mod4)$,\\
$|L|=\lceil\frac{n}{2}\rceil$ for $n\equiv3~(\mod4)$.

As above, if it is not the case of even $k$ and $n\equiv2~(\mod4)$, $(k-1)|L|$ is even. Then the graphs in $\mathcal{U}_{k,|L|}$ are $(k-1)$-regular. It follows that the graphs obtained from
 the complete bipartite graph with parts $L$ and $R$ by embedding a graph from
$\mathcal{U}_{k,|L|}$ in $G[L]$ and exactly one edge in $G[R]$ have the same equitable partition with three parts (see \cite{CDT2}). Thus they have the same spectral radius. Therefore, the authors \cite{CDT2} determine the exact spectral extremal graphs in this case. However, when $k$ is even and $n\equiv2~(\mod4)$, the structure of the graphs in ${\rm SPEX}(n,W_{2k+1})$ is not characterized completely. Our second result is to determine the exact spectral extremal graphs for this case.

For even $k\geq4$ and odd $n\geq k+1$, let $\mathcal{V}_{k,n}$ be the family of nearly $(k-1)$-regular graphs
of order $n$, in which one component is $K_{1}\vee\overline{M_{k-2}}\vee K_{2}$ and any other component has at most $2k-2$ vertices.

\begin{theorem} \label{main2}
Let $k\geq4$ be even and let $n\equiv2~(\mod4)$ be sufficiently large. If $G\in\rm{SPEX}(n,W_{2k+1})$, then $G$ is obtained from the Tur\'{a}n graph
with parts $L$ and $R$ by embedding a graph from
$\mathcal{V}_{k,\frac{n}{2}}$ in $G[L]$ and exactly one edge in $G[R]$.
\end{theorem}

Note that the graphs obtained from the Tur\'{a}n graph
with parts $L$ and $R$ by embedding a graph from
$\mathcal{V}_{k,\frac{n}{2}}$ in $G[L]$ and exactly one edge in $G[R]$, have the same equitable partition with six parts. Hence they have the same spectral radius. Therefore, Theorem \ref{main2} determines the exact spectral extremal graphs in the case.

The rest of this paper is organized as follows. In Section 2, we cite some lemmas from \cite{CDT2} and \cite{T}, and give the proof of Theorem \ref{main1}. In Section 3, we use a result in \cite{Z}, and give the proof of Theorem \ref{main2}.

\section{Proofs of Theorem \ref{main1}}

At the end of \cite{CDT2}, Cioab\u{a}, Desai and Tait pointed out that the main technical hurdle to prove Theorem \ref{CDT} for $k\in\left\{4,5\right\}$ is to prove Lemma 5.3 therein (i.e., the following Lemma \ref{main lemma}).
To prove Lemma \ref{main lemma}, we include some lemmas. The following one is taken from Theorem 4.4 of \cite{T}.

\begin{lem}{\rm (\cite{T})}\label{loop}
Let $H_{1}$ be a graph on $n_{0}$ vertices with maximum degree $d$ and $H_{2}$ be a
graph on $n-n_{0}$ vertices with maximum degree $d'$. $H_{1}$ and $H_{2}$ may have loops or multiple edges, where
loops add $1$ to the degree. Let $H$ be the join of $H_{1}$ and $H_{2}$. Define
\begin{center}B =
$\left(\begin{array}{cc}
d&n-n_{0}\\
n_{0}&d'
\end{array}\right)$.
   \end{center}
Then $\lambda_{1}(H)\leq\lambda_{1}(B)$.
\end{lem}

For $k\geq3$ and sufficiently large $n$, let $G\in {\rm SPEX}(n,W_{2k+1})$. Let $\mathbf{x}$ be the Perron vector of $G$ such that $\max_{u\in V(G)}x_{u}=1$. For two disjoint subsets $S,T$ of $V(G)$, let $e(S,T)=e_{G}(S,T)$.
The following Lemma \ref{max-cut} is taken from Lemma 3.2 and Lemma 3.6 of \cite{CDT2}.

\begin{lem}{\rm (\cite{CDT2})} \label{max-cut}
For any $\epsilon>0$, there is a partition $V(G)=S\cup T$ which forms a maximum cut satisfying $e(S,T)\geq(\frac{1}{4}-\epsilon)n^{2}$. Moreover,
$\frac{n}{2}-\sqrt{\frac{3nk}{2}}\leq |S|,|T|\leq\frac{n}{2}+\sqrt{\frac{3nk}{2}}$.
\end{lem}

The following Lemma \ref{degree} is taken from Lemma 3.5 and Lemma 5.1 of \cite{CDT2}.

\begin{lem}{\rm (\cite{CDT2})} \label{degree}
 $G[S]$ and $G[T]$ are $K_{1,k}$-free. Moreover, if $u$ is a vertex with degree $0\leq d\leq k-1$ in $G[S]$ (or $G[T]$), then $u$ is adjacent to all but at most $d$ vertices in $T$ (or $S$).
\end{lem}

The following Lemma \ref{entry} is taken from Lemma 3.7 of \cite{CDT2}.

\begin{lem}{\rm (\cite{CDT2})} \label{entry}
 For all $u\in V(G)$ and $0<\epsilon'<1$, we have $x_{u}>1-\epsilon'$.
\end{lem}

The following Fact 1 is taken from formula $(4)$ of \cite{CDT2}.

\medskip

\f{\bf Fact 1}. $\lambda_{1}(G)>\frac{k-1+\sqrt{(k-1)^{2}+n^{2}-1}}{2}+\frac{1}{2n}$.

\medskip

Then following Lemma \ref{max-degree} is taken from Lemma 5.2 of \cite{CDT2}.

\begin{lem}{\rm (\cite{CDT2})} \label{max-degree}
Assume that $C$ is a constant and that at most $C$ edges may be removed from $G$ so that the
graph induced by $S$ has maximum degree at most $a$ and the graph induced by $T$ has maximum degree at most $b$. Then
for $n$ large enough we must have $a+b\geq k-1$. If $C=0$ we must have $a+b\geq k$.
\end{lem}

\begin{lem}\label{main lemma}
For $k\geq3$, there exists at least one vertex in $G[S]$ or $G[T]$ with degree equal
to $k-1$.
\end{lem}

\f{\bf Proof:} The Claim 1 of Lemma 5.3 in \cite{CDT2} states that there exists at least one vertex in $G[S]$ or $G[T]$ with degree at least
 $\frac{k}{2}$. Then the lemma holds for $k=3$.  Let $L\in\left\{S,T\right\}$ be a part of $G$
such that $G[L]$ has a vertex of degree at least $\frac{k}{2}$. Without loss of generality, assume $L=S$. Let $R=T$. The Claim 2 of Lemma 5.3 in \cite{CDT2} states that there exists at least one vertex in $G[L]$  with degree at least $k-2$ (i.e., the following Claim 1). Using this result, the authors \cite{CDT2} confirmed the lemma for $k\geq6$.

\medskip

\f{\bf Claim 1}. (\cite{CDT2}) There exists at least one vertex  in $G[L]$ with degree at least $k-2$.

\medskip

Starting at Claim 1, we will prove the lemma for $k=4,5$.
By Claim 1, we can suppose that $$\max\left\{\Delta(G[L]),\Delta(G[R])\right\}=k-2,$$
 otherwise the lemma holds. We will obtain a contradiction in the following two cases.

\medskip

\f{\bf Case 1}. $k=4$. Then $\max\left\{\Delta(G[L]),\Delta(G[R])\right\}=k-2=2$, and thus each vertex in $G[L]$ and $G[R]$ has degree at most $2$.

\medskip

\f{\bf Claim 2}. If $\Delta(G[L])=2$ (or $\Delta(G[R])=2$), then  there are at most $11$ vertices of degree 2 in $G[R]$ (or $G[L]$).

\medskip

\f{\bf Proof of Claim 2}.
Assume that $v$ is a vertex in $G[L]$ with degree $2$. Let $v_{1}$ and $v_{2}$ be the two neighbors of $v$ in $L$. Then by Lemma \ref{degree}, $v,v_{1},v_{2}$ are all adjacent to a subset $R'\subseteq R$ of size $|R|-6$.
Since $G$ has no $W_{9}$, there are at most $5$ vertices contained in any two vertex-disjoint paths of $G[R']$. Thus there are no two disjoint paths of length  2 in $G[R']$. It follows that there are at most $5$ vertices of degree  2 in $G[R']$. Hence there are at most $11$ vertices of degree  2 in $G[R]$. Similarly, we can show that if there is a vertex of degree 2 in $G[R]$, then there are at most $11$ vertices of degree 2 in $G[L]$. This finishes the proof of Claim 2. \hfill$\Box$

Recall that $\max\left\{\Delta(G[L]),\Delta(G[R])\right\}=2$. By Claim 1 we have that $\Delta(G[L])=2$. If $\Delta(G[R])\leq1$, then we obtain a contradiction from Lemma \ref{max-degree} by letting $L=S,R=T$ and $a=2,b=1,C=0$. If $\Delta(G[R])=2$, by Claim 2, there are at most $22$ vertices of degree $2$ in $G[L]$ and $G[R]$. Thus, removing at most $22$ edges in $G[L]$ and $G[R]$ induce both of them to be graphs with maximum degree at most 1. Then we obtain a contradiction from Lemma \ref{max-degree} by letting $a=b=1$ and $C=22$.

\medskip

\f{\bf Case 2}. $k=5$.  Then $\max\left\{\Delta(G[L]),\Delta(G[R])\right\}=k-2=3$, and thus each vertex in $G[L]$ and $G[R]$ has degree at most $3$.

\medskip

\f{\bf Claim 3}. If $\Delta(G[L])=3$ and $\Delta(G[R])\geq2$ (or $\Delta(G[L])\geq2$ and $\Delta(G[R])=3$), then there are at most $22$ vertices of degree at least $2$ in $G[R]$ (or $G[L]$). Moreover, $G[L]$ (or $G[R]$) has at least $100$ components, in each of which there is at least one vertex with degree less than $3$.

\medskip

\f{\bf Proof of Claim 3}.
Firstly, assume that $\Delta(G[L])=3$ and $\Delta(G[R])\geq2$. Let $v$ be a vertex in $G[L]$ with degree $3$. Let $v_{1},v_{2}$ and $v_{3}$ be the three neighbors of $v$ in $L$. Then by Lemma \ref{degree}, $v,v_{1},v_{2},v_{3}$ are all adjacent to a subset $R'\subseteq R$ of size $|R|-12$.
Since $G$ has no $W_{11}$, there are at most $6$ vertices contained in any $3$ vertex-disjoint paths of $G[R']$. Thus there are no two disjoint paths of length  $2$ in $G[R']$. It follows that there are at most $10$ vertices of degree at least $2$ in $G[R']$. Hence there are at most $22$ vertices of degree at least  $2$ in $G[R]$.

Suppose that $u$ is a vertex in $G[R]$ with degree at least $2$. Let $u_{1}$ and $u_{2}$ be (any) two neighbors of $u$ in $R$. Then by Lemma \ref{degree}, $u,u_{1},u_{2}$ are all adjacent to a subset $L'\subseteq L$ of size $|L|-9$.
Since $G$ has no $W_{11}$, there are at most $7$ vertices contained in any two vertex-disjoint paths of $G[L']$. Hence there is no path of length $6$ in $G[L']$.  It follows that each component of $G[L']$ has diameter at most $5$. Recall that $\Delta(G[L'])\leq3$. Thus, each component of $G[L']$ has at most $94$ vertices.

Suppose that  $Q$ is a component of $G[L']$, in which each vertex has degree $3$. Then $|Q|\geq4$, and $Q$ has a path of order $4$. Since there are at most $7$ vertices contained in any $2$ vertex-disjoint paths of $G[L']$, we have that there is at most one such component $Q$  of $G[L']$, in which each vertex has degree $3$. Hence $G[L']$ has at least $\frac{|L|-9}{94}-1$ components $Q'$, in each of  which there is at least one vertex with degree less than $3$. Therefore, $G[L]$ has at least $\frac{|L|-9}{94}-1-27$ components $Q''$, in each of which there is at least one vertex with degree less than $3$. By Lemma \ref{max-cut} we see that $\frac{|L|-9}{94}-1-27\geq100$ for large $n$. Thus $G[L]$ has at least $100$ components, in each of which there is at least one vertex with degree less than $3$.

It is similar to prove for the case of $\Delta(G[L])\geq2$ and $\Delta(G[R])=3$. This finishes the proof of Claim 3. \hfill$\Box$

Recall that $\max\left\{\Delta(G[L]),\Delta(G[R])\right\}=3$. We have $\Delta(G[L])=3$ by Claim 1. If $\Delta(G[R])\leq1$, then we obtain  a contradiction from Lemma \ref{max-degree} by letting $a=3,b=1$ and $C=0$. If $\Delta(G[R])\geq2$, by Claim 3, there are at most $22$ vertices of degree at least 2 in $G[R]$. Moreover, $G[L]$ has such components $Q_{1},Q_{2},...,Q_{100}$, that for each $1\leq i\leq100$ there is at least one vertex, say $w_{i}$, with degree less than 3 in $Q_{i}$.  Since $G[R]$ has at most $22$ vertices of degree 2 or 3 by Claim 3,  there are at most $44$ edges of $G[R]$ whose deletion induces $G[R]$ to be a graph with maximum degree 1. Let $G'$ be the graph obtained from $G$ by deleting these $44$ edges, and adding an edge between $w_{2i-1}$ and $w_{2i}$ for any $1\leq i\leq 50$. Then $G'[L]$ has maximum degree 3 and $G'[R]$ has maximum degree 1. Let $\mathbf{x}^{T}$ be the transpose of $\mathbf{x}$. By Lemma \ref{entry}, for small $\epsilon'$ we have
$$\lambda_{1}(G')-\lambda_{1}(G)\geq \frac{\mathbf{x}^{T}A'(G)\mathbf{x}}{\mathbf{x}^{T}\mathbf{x}}-\frac{\mathbf{x}^{T}A(G)\mathbf{x}}{\mathbf{x}^{T}\mathbf{x}}
\geq\frac{2}{\mathbf{x}^{T}\mathbf{x}}(50(1-\epsilon')-44)>0,$$
implying that $\lambda_{1}(G')>\lambda_{1}(G)$. However, $G'$ is a spanning subgraph of the graph $G''$ obtained from $G'$ by adding all the non-edges between $L$ and $R$. Note that $\Delta(G''[L])=3$ and $\Delta(G''[R])=1$. Let
\begin{center}B =
$\left(\begin{array}{cc}
3&|R|\\
|L|&1
\end{array}\right)$.
   \end{center}
By Lemma \ref{loop}, we see that
$$\lambda_{1}(G'')\leq\lambda_{1}(B)=2+\sqrt{1+|L||R|}\leq2+\sqrt{1+\frac{n^{2}}{4}}.$$
 By Fact 1, we see $\lambda_{1}(G)>2+\sqrt{\frac{15}{4}+\frac{n^{2}}{4}}+\frac{1}{2n}$. However, $\lambda_{1}(G)<\lambda_{1}(G')\leq\lambda_{1}(G'')$. Then $2+\sqrt{\frac{15}{4}+\frac{n^{2}}{4}}+\frac{1}{2n}<2+\sqrt{1+\frac{n^{2}}{4}}$, a contradiction. This completes the proof . \hfill$\Box$

\medskip

\f{\bf Proof of Theorem \ref{main1}}.  Note that the conclusions from Lemma 5.4 to Lemma 5.8 of \cite{CDT2} hold for all $k\geq3$. Hence Theorem \ref{main1} follows. This completes the proof. \hfill$\Box$

 \section{Proof of Theorem \ref{main2}}

To prove Theorem \ref{main2}, we need to use a result in \cite{Z}. We first introduce some notations.
Let $G$ be a graph. For an integer $\ell\geq1$, a {\em walk} of length $\ell$ in $G$ is an ordered sequence of vertices $v_{0},v_{1},...,v_{\ell}$, such that $v_{i}\sim v_{i+1}$ for any $0\leq i\leq \ell-1$. The vertex $v_{0}$ is called the starting vertex of the walk. For any $u\in V(G)$, let $w^{\ell}_{G}(u)$ be the number of walks of length $\ell$ starting at $u$ in $G$. Let $W^{\ell}(G)=\sum_{v\in V(G)}w^{\ell}_{G}(u)$.
Clearly, $w^{1}_{G}(u)=d_{G}(u)$ for any $u\in V(G)$. For any integers $\ell\geq2$ and $1\leq i\leq \ell-1$, the following formula (by considering the $(i+1)$-th vertex in a walk of length $\ell$) will be used:
 $$W^{\ell}(G)=\sum_{u\in V(G)}w^{i}_{G}(u)w^{\ell-i}_{G}(u).$$

Let $G_{1}$ and $G_{2}$ be two graphs. We say $G_{1}\succ G_{2}$, if there is an integer $\ell\geq1$ such that $W^{\ell}(G_{1})>W^{\ell}(G_{2})$ and $W^{i}(G_{1})=W^{i}(G_{2})$ for any $1\leq i\leq \ell-1$;  $G_{1}\equiv G_{2}$, if $W^{i}(G_{1})=W^{i}(G_{2})$ for any $i\geq1$; $G_{1}\prec G_{2}$ , if $G_{2}\succ G_{1}$. 

Let $\mathcal{G}$ be a family of graphs.  
Let $${\rm EX}^{1}(\mathcal{G})=\left\{G\in\mathcal{G}~|~W^{1}(G)\geq W^{1}(G')~ {\rm for~any}~G'\in\mathcal{G}\right\},$$
 and  
$${\rm EX}^{\ell}(\mathcal{G})=\left\{G\in{\rm EX}^{\ell-1}(\mathcal{G})~|~W^{\ell}(G)\geq W^{\ell}(G')~ {\rm for~any}~G'\in{\rm EX}^{\ell-1}(\mathcal{G})\right\}$$
 for any $\ell\geq2$. Clearly, ${\rm EX}^{i+1}(\mathcal{G})\subseteq {\rm EX}^{i}(\mathcal{G})$ for any $i\geq1$. Let ${\rm EX}^{\infty}(\mathcal{G})=\cap_{1\leq i\leq\infty}{\rm EX}^{i}(\mathcal{G})$. 

\medskip

The following result is taken from \cite{Z}.

\begin{theorem}{\rm (Zhang \cite{Z})}\label{one-set}
Let $G$ be a connected graph of order $n$, and let $S$ be a subset of $V(G)$ with $1\leq|S|<n$. Assume that $T$ is a set of some isolated vertices of $G-S$, such that each vertex in $T$ is adjacent to each vertex in $S$ in $G$. Let $H_{1}$ and $H_{2}$ be two graphs with vertex set $T$. For any $1\leq i\leq2$, let $G_{i}$ be the graph obtained from $G$ by embedding the edges of $H_{i}$ into $T$. When $\lambda_{1}(G)$ is sufficiently large (compared with $|T|$), we have the following conclusions.\\
$(i)$ If $H_{1}\equiv H_{2}$, then $\lambda_{1}(G_{1})=\lambda_{1}(G_{2})$.\\
$(ii)$ If $H_{1}\succ H_{2}$, then $\lambda_{1}(G_{1})>\lambda_{1}(G_{2})$.\\
$(iii)$ If $H_{1}\prec H_{2}$, then $\lambda_{1}(G_{1})<\lambda_{1}(G_{2})$.
\end{theorem}

\begin{lem}\label{bounded order}
Let $G$ be a connected graph of order $n$, such that all vertices are of degree $\Delta$ except at most one vertex with degree $\Delta-1$, where $\Delta\geq2$. If $n\geq2\Delta+1$, then $G$ contains a path of order $2\Delta+1$.
\end{lem}

\f{\bf Proof:} We only prove the lemma for the case that there is a vertex with degree $\Delta-1$, since it is similar for the case that $G$ is $\Delta$-regular. Let $u$ be the unique vertex with degree $\Delta-1$ in $G$. The lemma holds trivially for $\Delta=2$. From now on, assume that $\Delta\geq3$. Let $P=v_{1},v_{2},...,v_{\ell}$ be a longest path in $G$. Then each neighbor of $v_{1}$ and $v_{\ell}$ is contained in $P$.
 If $\ell\geq2\Delta+1$, then the lemma holds. Thus we can assume that $\ell\leq2\Delta$. We will obtain a contradiction by two cases in the following.

\f{\bf Case 1.} $v_{1},v_{\ell}\neq u$.

Let $U=\left\{i-1~|~v_{1}\sim v_{i},2\leq i\leq \ell\right\}$ and $W=\left\{j~|~v_{\ell}\sim v_{j},1\leq j\leq \ell-1\right\}$. Then $|U|=|W|=\Delta$ as $d_{G}(v_{1})=d_{G}(v_{\ell})=\Delta$. Note that $U,W\subseteq\left\{1,2,...,\ell-1\right\}$ and $\ell\leq2\Delta$. Thus $|U\cap W|\geq|U|+|W|-(\ell-1)\geq1$.
Let $i_{0}\in U\cap W$. Then $v_{1}\sim v_{i_{0}+1}$ and  $v_{\ell}\sim v_{i_{0}}$.  It follows that $v_{1},v_{2},...,v_{i_{0}}v_{\ell},v_{\ell-1},...,v_{i_{0}+1}v_{1}$ is a cycle of order $\ell$. Since $G$ is a connected graph of order $n\geq2\Delta+1\geq\ell+1$, there is another vertex $v$ adjacent to some vertex in this cycle. Then we can obtain a path of order $\ell+1$. This is impossible as $P$ is a longest path in $G$.

\f{\bf Case 2.} Either $v_{1}$ or $v_{\ell} $ is $u$.

Without loss of generality, assume that $v_{1}=u$. Since $\Delta-1\geq2$ and all the neighbors of $v_{1}$ are contained in $P$, we have that $v_{1}$ is adjacent to $v_{j_{0}}$ for some $3\leq j_{0}\leq \ell$. Then $P'=v_{j_{0}-1},v_{j_{0}-2},...,v_{1},v_{j_{0}},v_{j_{0}+1},...,v_{\ell}$ is another longest path of order $\ell$ in $G$. Clearly, $v_{j_{0}-1},v_{\ell}\neq u$. Similar to the discussion in Case 1, we can obtain a contradiction.  This completes the proof. \hfill$\Box$

\medskip

By Lemma \ref{bounded order}, we see that $\mathcal{U}_{k,n}$ (defined in the introduction) is exactly the family of $(k-1)$-regular or nearly $(k-1)$-regular graphs
of order $n$ of which each component has at most $2k-2$ vertices.

For odd integers $\Delta\geq3$ and $n\geq3\Delta+4$, let $\mathcal{G}_{\Delta,n}$ be the set of graphs $G$ of order $n$, in which all vertices are of degree $\Delta$ except one vertex with degree $\Delta-1$, and all the components of $G$ are of order at most $2\Delta$. For even $k\geq4$ and odd $n\geq k+1$, recall that $\mathcal{V}_{k,n}$ is the family of nearly $(k-1)$-regular graphs
of order $n$, in which one component is $K_{1}\vee\overline{M_{k-2}}\vee K_{2}$ and any other component has at most $2k-2$ vertices. Clearly, $\mathcal{V}_{k,n}\subseteq \mathcal{G}_{k-1,n}= \mathcal{U}_{k,n}$.

\begin{lem}\label{walk lemma}
Let $\Delta\geq3$ and $n\geq3\Delta+4$ be odd integers. Then ${\rm EX}^{\infty}(\mathcal{G}_{\Delta,n})=\mathcal{V}_{\Delta+1,n}$.
\end{lem}

\f{\bf Proof:} Note that there exists a $\Delta$-regular graph of order $m$ if and only if $\Delta m$ is even and $m\geq\Delta+1$. Since $|K_{1}\vee\overline{M_{\Delta-1}}\vee K_{2}|=\Delta+2$ and $n-(\Delta+2)\geq2(\Delta+1)$ can be partitioned into several even integers between $\Delta+1$ and $2\Delta$, we see that $\mathcal{V}_{\Delta+1,n}$
is not empty. Let $H\in\mathcal{V}_{\Delta+1,n}$. Then the component of $H$ including the vertex with degree $\Delta-1$ is $K_{1}\vee\overline{M_{\Delta-1}}\vee K_{2}$.

Let $G\in {\rm EX}^{\infty}(\mathcal{G}_{\Delta,n})$. Let $u$ be the vertex with degree $\Delta-1$ in $G$, and let $Q$ be the component of $G$ including $u$. We will prove $Q=K_{1}\vee\overline{M_{\Delta-1}}\vee K_{2}$. Then $G\in\mathcal{V}_{\Delta+1,n}$, and the lemma follows. (Note that all the graphs in $\mathcal{V}_{\Delta+1,n}$ have the same number of walks of length $\ell$ for any $\ell\geq1$.) 

Set $q=|Q|$. Then $Q$ has one vertex with degree $\Delta-1$ and $q-1$ vertices with degree $\Delta$. Note that $q\leq2\Delta$ is odd. Then $q\leq2\Delta-1$ and $n-q\geq\Delta+5$ is even. If the diameter of $Q$ is at least 3, then $q\geq\Delta+\Delta+1$, a contradiction. Thus the diameter of $Q$ is 2. For any $1\leq i\leq 2$, let $N_{i}$ be the set of vertices of $Q$ at distance $i$ from the vertex $u$. Then $|N_{1}|=\Delta-1$ and $|N_{2}|=q-(1+\Delta-1)\leq\Delta-1$. Thus $2\leq |N_{2}|\leq\Delta-1$ as $|N_{2}|$ is even. Let $e(N_{1},N_{2})=e_{G}(N_{1},N_{2})$.
It follows that
 $$e(N_{1},N_{2})\geq\sum_{v\in N_{2}}(\Delta+1-|N_{2}|)|N_{2}|\geq2(\Delta-1),$$
  with equality if and only if $Q[N_{2}]$ is a complete graph of order $2$ or $\Delta-1$. 
  
  For any $v\in V(Q)$ and $1\leq i\leq2$, let $d_{i}(v)$ be the number of neighbors of $v$ in $N_{i}$. Let $w^{i}(u)=w^{i}_{G}(u)$ for any $i\geq1$ and $u\in V(G)$. It is easy to check that:
 $$w^{2}(u)=\Delta^{2}-\Delta,w^{3}(u)=(\Delta^{2}-1)(\Delta-1);$$
for $v\in N_{1}$,
$$w^{2}(v)=\Delta^{2}-1,w^{3}(v)=\Delta^{2}-\Delta+d_{1}(v)(\Delta^{2}-1)+d_{2}(u)\Delta^{2}=\Delta^{3}-2\Delta+1+d_{2}(v);$$
 for $v\in N_{2}$,
$$w^{2}(v)=\Delta^{2},w^{3}(v)=d_{1}(v)(\Delta^{2}-1)+d_{2}(u)\Delta^{2}=\Delta^{3}-d_{1}(v);$$
for $v\in V(G)-V(Q)$ and $i\geq1$,
$$w^{i}(v)=\Delta^{i}.$$
Then by a calculation, we have that
$$W^{1}(G)=n\Delta-1,$$
$$W^{2}(G)=(n-1)\Delta^{2}+(\Delta-1)^{2},$$
$$W^{3}(G)=\sum_{v\in V(G)}w^{1}(v)w^{2}(v)=n\Delta^{3}-3\Delta^{2}+2\Delta,$$
$$W^{4}(G)=\sum_{v\in V(G)}w^{2}(v)w^{2}(v)=n\Delta^{4}-4\Delta^{3}+3\Delta^{2}+\Delta-1,$$
$$W^{5}(G)=\sum_{v\in V(G)}w^{2}(v)w^{3}(v)=(n-q)\Delta^{5}+\sum_{v\in V(Q)}w^{2}(v)w^{3}(v),$$
$$=(n-\Delta)\Delta^{5}+\Delta(\Delta+1)(\Delta-1)^{3}+(\Delta^{2}-1)(\Delta-1)(\Delta^{3}-2\Delta+1)-e(N_{1},N_{2})$$
and
$$W^{6}(G)=\sum_{v\in V(G)}w^{3}(v)w^{3}(v)=(n-q)\Delta^{6}+\sum_{v\in V(Q)}w^{3}(v)w^{3}(v)$$
$$=(n-\Delta)\Delta^{6}+(\Delta-1)^{2}(\Delta^{2}-1)^{2}+(\Delta-1)(\Delta^{3}-2\Delta+1)^{2}-(4\Delta-2)e(N_{1},N_{2})
$$
$$+\sum_{v\in N_{1}}d_{2}^{2}(v)+\sum_{v\in N_{2}}d_{1}^{2}(v).$$

Clearly, $W^{i}(G)=W^{i}(H)$ is constant for $1\leq i\leq4$. Since $G\in {\rm EX}^{\infty}(\mathcal{G}_{\Delta,n})$, we have $W^{5}(G)\geq W^{5}(H)$. This requires that $e(N_{1},N_{2})\geq2(\Delta-1)$ is minimized. Then $e(N_{1},N_{2})=2(\Delta-1)$ and $W^{5}(G)= W^{5}(H)$. Moreover, $G[N_{2}]$ is a complete graph of order $2$ or $\Delta-1$. Also by $G\in {\rm EX}^{\infty}(\mathcal{G}_{\Delta,n})$, we have $W^{6}(G)\geq W^{6}(H)$ as $W^{i}(G)= W^{i}(H)$ for $1\leq i\leq5$.
This requires that $\sum_{v\in N_{1}}d_{2}^{2}(v)+\sum_{v\in N_{2}}d_{1}^{2}(v)$ is maximized. 

Now we show that $\sum_{v\in N_{1}}d_{2}^{2}(v)+\sum_{v\in N_{2}}d_{1}^{2}(v)$ is maximized when $G[N_{2}]=K_{2}$. Recall that $e(N_{1},N_{2})=2(\Delta-1)$, and $G[N_{2}]=K_{2}$ or $G[N_{2}]=K_{\Delta-1}$. If $\Delta=3$, then we are done. Now suppose that $\Delta\geq5$. If $G[N_{2}]=K_{2}$, then each vertex in $N_{1}$ has at most two neighbors in $N_{2}$. Since $e(N_{1},N_{2})=2(\Delta-1)$ and $|N_{1}|=\Delta-1$, we see that each vertex in $N_{1}$ has exactly two neighbors in $N_{2}$, and thus $Q=K_{1}\vee\overline{M_{\Delta-1}}\vee K_{2}$. It follows that $$\sum_{v\in N_{1}}d_{2}^{2}(v)+\sum_{v\in N_{2}}d_{1}^{2}(v)=2\Delta^{2}-2.$$
 If $G[N_{2}]=K_{\Delta-1}$, then $\sum_{v\in N_{2}}d_{1}^{2}(v)=4\Delta-4$.
Since $1\leq d_{2}(v)\leq\Delta-1$ for each $v\in N_{1}$ and $\sum_{v\in N_{1}}d_{2}(v)=2(\Delta-1)$, we have 
$$\sum_{v\in N_{1}}d_{2}^{2}(v)<(\Delta-1)^{2}+(\Delta-1)^{2}=2\Delta^{2}-4\Delta+2.$$
It follows that 
$$\sum_{v\in N_{1}}d_{2}^{2}(v)+\sum_{v\in N_{2}}d_{1}^{2}(v)<4\Delta-4+2\Delta^{2}-4\Delta+2=2\Delta^{2}-2.$$
Thus $\sum_{v\in N_{1}}d_{2}^{2}(v)+\sum_{v\in N_{2}}d_{1}^{2}(v)$ is maximized when $G[N_{2}]=K_{2}$.
Consequently, we must have $Q=K_{1}\vee\overline{M_{\Delta-1}}\vee K_{2}$. This completes the proof. \hfill$\Box$

\medskip

\f{\bf Proof of Theorem \ref{main2}}. Let $k\geq4$ be even and let $n\equiv2~(\mod4)$ be sufficiently large. Assume that $G\in\rm{SPEX}(n,W_{2k+1})$. By Theorem \ref{CDT}, $G$ is obtained from
 a complete bipartite graph with parts $L$ and $R$ of size $\frac{n}{2}+s$ and $\frac{n}{2}-s$ respectively, by embedding a graph from
$\mathcal{U}_{k,\frac{n}{2}+s}$ in $G[L]$ and exactly one edge in $G[R]$. Moreover, $|L|\in\left\{\frac{n}{2},\frac{n}{2}+1\right\}$.

\f{\bf Claim 1}. $|L|=\frac{n}{2}$. 

\f{\bf Proof of  Claim 1}. Let $H$ be any graph obtained from the Tur\'{a}n graph
with parts $L$ and $R$ (of equal size $\frac{n}{2}$) by embedding a graph from
$\mathcal{V}_{k,\frac{n}{2}}$ in $G[L]$ and exactly one edge in $G[R]$. Then the component of $H[L]$ including the vertex with degree $k-2$ is $K_{1}\vee\overline{M_{k-2}}\vee K_{2}$. Clearly, $V(H)$ has an equitable partition of $6$ parts: $$V(K_{1}),V(\overline{M_{k-2}}),V(K_{2}),L-V(K_{1})\cup V(\overline{M_{k-2}})\cup V(K_{2}),V(K_{2}),R-V(K_{2}).$$
The quotient matrix of this partition is
\begin{center}$B_{1}=$
$\left(\begin{array}{cccccc}
 0&k-2&0&0&2&\frac{n}{2}-2\\
1&k-3&2&0&2&\frac{n}{2}-2\\
0&k-2&1&0&2&\frac{n}{2}-2\\
0&0&0&k-1&2&\frac{n}{2}-2\\
1&k-2&2&\frac{n}{2}-k-1&1&0\\
1&k-2&2&\frac{n}{2}-k-1&0&0
\end{array}\right)$.
   \end{center}
Then $\lambda_{1}(H)=\lambda_{1}(B_{1})$.
Let $H'$ be any graph obtained from a complete bipartite graph with parts $L$ and $R$ of size $\frac{n}{2}+1$ and $\frac{n}{2}-1$ respectively, by embedding a graph from
$\mathcal{U}_{k,\frac{n}{2}+1}$ in $G[L]$ and exactly one edge in $G[R]$. Clearly, $V(H')$ has an equitable partition of $3$ parts: $L, V(K_{2}),R-V(K_{2})$. The quotient matrix of this partition is
\begin{center}$B_{2}=$
$\left(\begin{array}{ccc}
k-1&2&\frac{n}{2}-3\\
\frac{n}{2}+1&1&0\\
\frac{n}{2}+1&0&0
\end{array}\right)$.
   \end{center}
Then $\lambda_{1}(H')=\lambda_{1}(B_{2})$. Let $f(B_{i},x)$ be the characteristic polynomial of $B_{i}$ for $i=1,2$. By a calculation, we have
$$f(B_{1},x)-(x^{3}+(2-k)x^{2}+(2-2k)x+\frac{3-k}{2}n+k-1)f(B_{2},x)=n^{4}(-\frac{k-3}{8}\frac{x}{n}-\frac{1}{4}(\frac{x}{n})^{2}+\Theta(\frac{1}{n})).$$
Note that $f(B_{2},\lambda_{1}(B_{2}))=0$. It follows that $f(B_{1},\lambda_{1}(B_{2}))<0$ for large $n$. This implies that $\lambda_{1}(B_{1})>\lambda_{1}(B_{2})$, i.e., $\lambda_{1}(H)>\lambda_{1}(H')$. Since $G\in\rm{SPEX}(n,W_{2k+1})$, we must have $|L|=\frac{n}{2}$. This finishes the proof of Claim 1. \hfill$\Box$

By Claim 1, we see that $G$ is obtained from
the Tur\'{a}n graph with parts $L$ and $R$, by embedding a graph from
$\mathcal{U}_{k,\frac{n}{2}}$ in $G[L]$ and exactly one edge in $G[R]$. Let $Q_{1}$ be the component of $G[L]$ including the vertex with degree $k-2$.
Since each component of $G[L]$ has at most $2k-2$ vertices and at least $k$ vertices, we can choose several components $Q_{2},..,Q_{t}$, such that $3k+1\leq|\cup_{1\leq i\leq t}Q_{i}|\leq5k$. Let $Q=\cup_{1\leq i\leq t}Q_{i}$. Since $G\in\rm{SPEX}(n,W_{2k+1})$, by Theorem \ref{one-set} we must have $Q\in{\rm EX}^{\infty}(\mathcal{G}_{k-1,|Q|})$ as $n$ is large with respective to $|Q|$. By Lemma \ref{walk lemma} we have $Q\in\mathcal{V}_{k,|Q|}$. Hence $G$ is obtained from the Tur\'{a}n graph
with parts $L$ and $R$ by embedding a graph from
$\mathcal{V}_{k,\frac{n}{2}}$ in $G[L]$ and exactly one edge in $G[R]$.  This completes the proof. \hfill$\Box$

\medskip

\f{\bf Declaration of competing interest}

\medskip

There is no conflict of interest.

\medskip

\f{\bf Data availability statement}

\medskip

No data was used for the research described in the article.


\medskip

\begin{thebibliography}{99}









\bibitem{BDT}
 J. Byrne, D.N. Desai and M. Tait, A general theorem in spectral extremal graph theory,
arXiv:2401.07266v1.

\bibitem{CRS}
D. Cvetkovi\'{c}, P. Rowlinson and S. Simi\'{c}, An Introduction to the Theory of Graph Spectra, Cambridge University Press, Cambridge, 2010.

\bibitem{CFTZ}
S. Cioab\u{a}, L. Feng, M. Tait and X. Zhang, The maximum spectral radius of graphs
without friendship subgraphs, Electron. J. Combin. 27 (2020)  \#4.22.

\bibitem{CDT1}
S. Cioab\u{a}, D. Desai and M. Tait, The spectral even cycle problem, arXiv:2205.00990v1.


\bibitem{CDT2}
S. Cioab\u{a}, D. Desai and M. Tait, The spectral radius of graphs with no odd wheels,
European J. Combin. 99 (2022) 103420.

\bibitem{CDT3}
S. Cioab\u{a}, D. Desai and M. Tait,
A spectral Erd\H{o}s-S\'{o}s theorem, SIAM J. Discrete
Math. 37 (2023) 2228-2239.

\bibitem{DKLNTW}
 D. Desai, L. Kang, Y. Li, Z. Ni, M. Tait and J. Wang, Spectral extremal graphs
for intersecting cliques, Linear Algebra Appl. 644 (2022) 234-258.

\bibitem{DJ}
 T. Dzido and A. Jastrzebski, Tur\'{a}n numbers for odd wheels, Discrete Math. 341 (4) (2018) 1150-1154.

\bibitem{FLSZ}
L. Fang, H. Lin, J. Shu and Z. Zhang, Spectral extremal results on trees, arXiv:2401.05786v2.
\bibitem{FTZ}
L. Fang, M. Tait and M. Zhai, Tur\'{a}n numbers for non-bipartite graphs and applications to spectral extremal problems, arXiv:2404.09069v1.

\bibitem{LL}
 X. Lei and S. Li, Spectral extremal problem on disjoint color-critical graphs, Electron. J. Combin. 31 (2024) \# 1.25.

\bibitem{LP1}
Y. Li, and Y. Peng, Refinement on spectral Tur\'{a}n's theorem, SIAM J. Discrete Math.
37 (2023) 2462-2485.
\bibitem{LP2}
Y. Li, and Y. Peng, The spectral radius of graphs with no intersecting odd cycles,
Discrete Math. 345 (2022) 112907.

\bibitem{LZZ}
 H. Lin, M. Zhai and Y. Zhao, Spectral radius, edge-disjoint cycles and cycles of
the same length, Electron. J. Combin. 29 (2022) \# 2.1.

\bibitem{NWK}
 Z. Ni, J. Wang and L. Kang, Spectral extremal graphs for disjoint cliques, Electron.
J. Combin. 30 (2023)  \#1.20.

\bibitem{N1}
V. Nikiforov, The spectral radius of graphs without paths and cycles of specified length, Linear Algebra Appl. 432 (2010) 2243-2256.

\bibitem{N2}
 V. Nikiforov, Spectral saturation: inverting the spectral Tur\'{a}n theorem, Electron. J. Combin. 16 (R33) (2009) 1.

\bibitem{S}
 M. Simonovits, A method for solving extremal problems in graph theory, stability problems, in: Theory of Graphs
(Proc. Colloq., Tihany, 1966), 1968, pp. 279-319.


\bibitem{T}
M. Tait, The Colin de Verdi\`{e}re parameter, excluded minors, and the spectral radius, J. Combin. Theory Ser. A 166 (2019) 42-58.

\bibitem{WKX}
J. Wang, L. Kang and Y. Xue, On a conjecture of spectral extremal problems, J.
Combin. Theory Ser. B 159 (2023) 20-41.

\bibitem{Y}
L. Yuan, Extremal graphs for wheels, J.
Graph Theory 98 (2021) 691-707.

\bibitem{ZHL}
Y. Zhao, X. Huang and H. Lin, The maximum spectral radius of wheel-free graphs, Discrete Math. 344 (5)
(2021) 112341.

\bibitem{ZL}
 M. Zhai and H. Lin, Spectral extrema of graphs: Forbidden hexagon, Discrete Math.
343 (2020) 112028.

\bibitem{Z}
W. Zhang, Walks, infinite series and spectral radius of graphs, arXiv:2406.07821v1.
\end{thebibliography}
\end{document}